\documentstyle[12pt]{article}
\newcommand{\sect}[1]{\section {#1}\setcounter{equation}{0}}

\font\mbn=msbm10 scaled \magstep1
\font\mbs=msbm7 scaled \magstep1
\font\mbss=msbm5 scaled \magstep1
\newfam\mbff
\textfont\mbff=\mbn
\scriptfont\mbff=\mbs
\scriptscriptfont\mbff=\mbss
\def\mbf{\fam\mbff}\def\Re{{\mbf R}}
\def\Co{{\mbf C}}

\def\Z{{\mbf Z}}

\def\To{{\mbf T}}
\newtheorem{Th}{Theorem}[section]
\newtheorem{Lm}[Th]{Lemma}

\newtheorem{Prp}[Th]{Proposition}
\newtheorem{R}[Th]{Remark}
\author{Alexander Brudnyi\thanks
{Research supported in part by NSERC.
}}
\title{HOLOMORPHIC FUNCTIONS OF EXPONENTIAL GROWTH ON ABELIAN COVERINGS OF
A PROJECTIVE MANIFOLD}
\date{}
\begin{document}
\maketitle
%=====================================================
\begin{abstract}
Let $M$ be a projective manifold, $p:M_{G}\longrightarrow M$ a
regular covering over $M$ with a free abelian transformation group
$G$. We describe holomorphic functions on $M_{G}$ of an
exponential growth with respect to the distance defined by a
metric pulled back from $M$. As a corollary we obtain for such
functions Cartwright and Liouville type theorems. Our approach
brings together $L_{2}$ cohomology technique for holomorphic
vector bundles on complete K\"{a}hler manifolds and geometric
properties of projective manifolds.
\end{abstract}
%==============================================================
%=================================================
\sect{\hspace*{-1em}. Introduction.} {\bf 1.1.} Recently there was
an essential progress in study of harmonic functions of polynomial
growth on complete Riemannian manifolds (see, in particular, [CM],
[Gu], [Ka], [L], [LZ], [Li], [LySu] for the results and further
references). As a corollary one also obtains a description of
holomorphic functions of polynomial growth on nilpotent coverings
of compact K\"{a}hler manifolds (see also [Br]). On the other
hand, very little is known about existence and behaviour of slowly
growing harmonic (respectively holomorphic) functions on covering
spaces of compact Riemannian (respectively K\"{a}hler) manifolds.
The methods of the above cited papers seem to be not sufficient
for application to the general situation. This paper is devoted to
study of slowly growing holomorphic functions on abelian coverings
of projective manifolds. Our approach is based on $L_{2}$
cohomology technique for holomorphic vector bundles on complete
K\"{a}hler manifolds and geometric properties of projective
manifolds and differs from the methods of the above mentioned
papers.

In order to formulate the results of the paper we consider a
projective manifold $M$ and its regular covering
$p:M_{G}\longrightarrow M$ with a free abelian transformation
group $G$. Denote by $r$ the distance from a fixed point in
$M_{G}$ defined by a metric pulled back from $M$. We study
holomorphic functions $f$ on $M_{G}$ satisfying (for some
$\epsilon>0$)
\begin{equation}\label{1}
|f(z)|\leq ce^{\epsilon r^{2}(z)},\ \ \ \ (z\in M_{G}).
\end{equation}
Recall that the covering space $M_{G}$ can be described as
follows.

Let $\omega_{1},...,\omega_{n}$ be a basis of holomorphic 1-forms
on $M$ and $A:M\longrightarrow\Co\To^{n}$ be the Albanese map of
$M$ associated with this basis. By definition, $$
A(z)=\left(\int_{z_{0}}^{z}\omega_{1},...,\int_{z_{0}}^{z}\omega_{n}\right)
$$ for a fixed $z_{0}\in M$. Consider a free abelian  quotient
group $G$ of the fundamental group
$\pi_{1}(\Co\To^{n})\cong\Z^{2n}$. Let
$t:T_{G}\longrightarrow\Co\To^{n}$ be the regular covering over
torus with the transformation group $G$. We can think of $T_{G}$
as a locally trivial fibre bundle over $\Co\To^{n}$ with discrete
fibres. Then $M_{G}=A^{*}T_{G}$ is the pullback of $T_{G}$ to $M$.
By definition the fundamental group of $M_{G}$ is $H:=(\pi\circ
A_{*})^{-1}(G)\subset\pi_{1}(M)$, where
$\pi:\Z^{2n}\longrightarrow G$ denotes the quotient map. By the
covering homotopy theorem there is a proper holomorphic map
$A_{G}:M_{G}\longrightarrow T_{G}$ that covers $A$ and such that
$\widetilde M_{G}:=A_{G}(M_{G})\subset T_{G}$ is a covering of
complex variety $A(M)\subset\Co\To^{n}$.

Our main result shows that if $f$ satisfies (\ref{1}) then there
is a uniquely defined holomorphic function $g$ on $T_{G}$ with a
similar growth condition such that $f=A_{G}^{*}(g)$. To its
formulation we let $\phi$ be a smooth nonnegative function on
$T_{G}$ and $\tilde\phi=A_{G}^{*}(\phi)$. Consider the Hilbert
space ${\cal H}_{\tilde\phi}(M_{G})$ of holomorphic functions $f$
on $M_{G}$ with the norm $$
|f|:=\int_{M_{G}}|f|^{2}e^{-\tilde\phi}dV. $$ Here $dV$ is the
pullback of the volume form on $M$ defined by a K\"{a}hler metric.
Similarly we introduce the Hilbert space ${\cal H}_{\phi}(T_{G})$
of holomorphic functions $f$ on $T_{G}$ with the norm $$
|f|:=\int_{T_{G}}|f|^{2}e^{-\phi}d\widetilde V, $$ where
$d\widetilde V$ is the pullback of the standard volume form on
$\Co\To^{n}$. Let $\{dz_{1},...,dz_{n}\}$ be the basis of
holomorphic 1-forms on $\Co\To^{n}$ such that
$A^{*}(dz_{i})=\omega_{i}$ for $i=1,...,n$. By the same symbol we
denote the pullback of these forms to $T_{G}$. Let ${\cal
L}(\phi)= \sum_{i,j}a_{ij}(z,\overline{z})dz_{i}\wedge
d\overline{z}_{j}$ be the Levi form of $\phi$. We set $$ |{\cal
L}(\phi)|:=\sup_{i,j,z\in T_{G}}|a_{ij}(z)|\ . $$

Assume that there is a constant $c>0$ such that
\begin{equation}\label{lip}
|\phi(x)-\phi(y)|\leq cd(x,y),
\end{equation}
where $d(.,.)$ is the distance on $T_{G}$ defined by the pullback of
the flat metric on $\Co\To^{n}$.
\begin{Th}\label{te1}
There is a constant $C=C(M,A)>0$ such that if $|{\cal L}(\phi)|<C$
and $\phi$ satisfies (\ref{lip}) then $A_{G}^{*}$ maps ${\cal
H}_{\phi}(T_{G})$ isomorphically onto ${\cal
H}_{\tilde\phi}(M_{G})$.
\end{Th}

Assume now that instead of (\ref{lip}) $\phi$ satisfies:

for any $\epsilon>0$, $x,y\in T_{G}$ with $d(x,y)\leq t$ there is
a function $c(\epsilon,t)>0$ increasing in $t$ such that
\begin{equation}\label{ineq}
\phi(x)\leq (1+\epsilon)\phi(y)+c(\epsilon,t)\ .
\end{equation}
\begin{Th}\label{te2}
Let $C$ be as in Theorem \ref{te1}, $|{\cal L}(\phi)|<C'<C$ and
$\phi$ satisfies (\ref{ineq}). There is a constant $\tilde\epsilon
(C')>0$ such that for any $f\in {\cal H}_{\tilde\phi}(M_{G})$
there exists a unique $\hat f\in\cap_{\epsilon<\tilde\epsilon
(C')} {\cal H}_{(1+\epsilon)\phi}(T_{G})$ satisfying $$
A_{G}^{*}(\hat f)=f\ \ \ {\rm and}\ \ \ |\hat f|\leq
c(\epsilon)|f|.$$  Here we regard $\hat f$ as an element of ${\cal
H}_{(1+\epsilon)\phi}(T_{G})$.
\end{Th}

In the following examples $M_{G}$ is a regular covering over $M$
with the maximal free abelian transformation group $G$ (so
$T_{G}=\Co^{n}$).\\ {\bf Examples.} 1. Let
$\phi(z)=k\log(p+|z|^{2})$ on $\Co^{n}$, where $|z|$ is the
Euclidean norm of the vector $z\in\Co^{n}$ and $p>0$ is so big
that $|{\cal L}(\phi)|<C$. (Such $p$ exists because ${\cal L}(\log
|z|)\to 0$ when $|z|\to\infty$.) Then ${\cal H}_{\phi}(\Co^{n})$
is isomorphic to the space of holomorphic polynomials of degree
$\leq k-n-1$. Therefore every holomorphic function on $M_{G}$ of
the corresponding polynomial growth is the pullback by $A_{G}$ of
a uniquely defined holomorphic polynomial on $\Co^{n}$. This gives
another proof for projective manifolds of the main result of [Br].
\\
2. Let $\phi(z)=2\sigma\sqrt{p+|z|^{2}}$ on
$\Co^{n}$, where $p$ is such that $|{\cal L}(\phi)|<C$. Then
${\cal H}_{\phi}(\Co^{n})$ consists of entire functions of the
exponential type $<\sigma$. Now Theorem \ref{te1} describes
holomorphic functions $f$ on $M_{G}$ satisfying
$|f(z)|<ce^{\sigma' r(z)}$, $z\in M_{G}$, $\sigma'<\sigma$.\\ 3.
Let $\phi(z)=2\sigma|z|^{2}$ on $\Co^{n}$ with $2\sigma<C$. Then the
assumptions of
Theorem \ref{te2} are fulfilled and the theorem describes holomorphic
functions $f$ on $M_{G}$ satisfying $|f(z)|<ce^{\sigma' r^{2}(z)}$, $z\in
M_{G}$, $\sigma'<\sigma$.\\
4. Assume that $C$ is a compact complex curve of genus $g\geq 1$. Then
$C_{G}$  can be thought of as a submanifold in $\Co^{g}$.
Applying Theorem \ref{te2} we obtain the following Cartwright type theorem.

{\em There is a positive number $\sigma=\sigma (C_{G})$ such that
any holomorphic function $f$ on $\Co^{g}$ satisfying $|f(z)|\leq
ce^{\sigma'|z|^{2}}$, $0<\sigma'<\sigma$, $z\in\Co^{g}$, and
$f|_{C_{G}}=0$, equals 0 identically.}\\
\\
{\bf 1.2.} The classical Liouville theorem asserts that every
bounded holomorphic function on $\Co^{n}$ is a constant. Based on
Theorem \ref{te1} we prove Liouville type theorems for holomorphic
functions of slow growth on abelian coverings over a projective
manifold.

Let $\Gamma\subset
H_{1}(M,\Z)\cong\pi_{1}(M)/[\pi_{1}(M),\pi_{1}(M)]$ be the maximal
free abelian subgroup of the homology group of $M$. Further, let
$\Omega^{1}(M)$ be the space of holomorphic 1-forms on $M$. Any
$\omega\in\Omega^{1}(M)$ determines a complex-valued linear
functional on $\Gamma$ by integration. For a subgroup
$H\subset\Gamma$ denote by $\Lambda (H)$ the minimal complex
subspace of holomorphic 1-forms vanishing on $H$. Assume also that
the quotient group $G=\Gamma/H$ is torsion free and $M_{G}$ is the
regular covering over $M$ with the transformation group $G$.
\begin{Th}\label{te3}
Let $H$ be such that $\Lambda (H)=\Omega^{1}(M)$. Then any
holomorphic on $M_{G}$ function $f$ satisfying for any
$\epsilon>0$ $$ |f(z)|\leq c(\epsilon)e^{\epsilon r(z)}\ \ \ \
(z\in M_{G}) $$ is a constant.
\end{Th}
\begin{R}\label{con}
It can be conjectured that the results of this paper are also true
for abelian coverings of an arbitrary compact K\"{a}hler manifold.
\end{R}
%====================================
\sect{\hspace*{-1em}. Preliminaries.} {\bf 2.1. $L_{2}$ cohomology
theory.} In the proof of our main results we use $L_{2}$
cohomology technique for holomorphic vector bundles on complete
K\"{a}hler manifolds. We start by reviewing some results of
$L_{2}$ cohomology (see, e.g., L\'{a}russon [La] for more details
and further references).

Let $X$ be a complex manifold of dimension $n$ with a hermitian
metric and $E$ be a holomorphic vector bundle over $X$ with a
hermitian metric. Let $L_{2}^{p,q}(X,E)$ be the space of
$E$-valued $(p,q)$-forms on $X$ with the $L_{2}$ norm, and let
$W_{2}^{p,q}(X,E)$ be the subspace of forms $\eta$ such that
$\overline{\partial}\eta$ is $L_{2}$. The forms $\eta$ may be
taken to be either smooth or just measurable, in which case
$\overline{\partial}\eta$ is understood in the distributional
sense. The cohomology of the resulting $L_{2}$ Dolbeault complex
$(W_{2}^{\cdot,\cdot}, \overline{\partial})$ is the
$L_{2}$-cohomology $$
H_{(2)}^{p,q}(X,E)=Z_{2}^{p,q}(X,E)/B_{2}^{p,q}(X,E), $$ where
$Z_{2}^{p,q}(X,E)$ and $B_{2}^{p,q}(X,E)$ are the spaces of
$\overline{\partial}$-closed and $\overline{\partial}$-exact forms
in $L_{2}^{p,q}(X,E)$, respectively. Let $E^{*}$ be the dual
bundle of $E$ with the dual metric. In our proofs we use the
following result discovered by L\'{a}russon [La].
\begin{Prp}\label{lar}
Let $E$ be a hermitian vector bundle with curvature $\Theta$ on a complex
manifold $X$ of dimension $n\geq 2$ with a complete K\"{a}hler form $\omega$.
If $\Theta\geq\epsilon\omega$ for some $\epsilon>0$ in the sense of
Nakano, then
$$
H_{(2)}^{0,q}(X,E^{*})=0\ \ \ \ \ {\rm for}\ q<n\ .
$$
\end{Prp}
\begin{R}\label{lin}
Let $E$ satisfy conditions of Proposition \ref{lar}. Consider
linear map $\overline{\partial}:W_{2}^{0,0}(X,E^{*})\longrightarrow
Z_{2}^{0,1}(X,E^{*})$ and introduce the norm in $W_{2}^{0,0}(X,E^{*})$ by
$$
|f|:=|f|_{2}+|\overline{\partial}f|_{2},\ \ \ \ \ f\in W_{2}^{0,0}(X,E^{*})\ .
$$
According to Proposition \ref{lar} for $q=1$ and $q=0$, there is a linear
map $s:Z_{2}^{0,1}(X,E^{*})\longrightarrow W_{2}^{0,0}(X,E^{*})$ such that
$s\circ\overline{\partial}=id$ and $\overline{\partial}\circ s=id$.
Then by the Banach theorem, $\overline{\partial}$ is open and
$s=(\overline{\partial})^{-1}$.
\end{R}
{\bf 2.2. $\overline{\partial}$-method.}
Let $i:X\hookrightarrow Y$ be a complex compact submanifold of
codimension 1 of an $n$-dimensional compact K\"{a}hler manifold $Y$,
$n\geq 2$, with a K\"{a}hler form $\omega$. Assume that the induced
homomorphism $i_{*}:H_{1}(X,\Re)\longrightarrow H_{1}(Y,\Re)$ is surjective.
Let $G$ be a free abelian quotient group of $\pi_{1}(Y)$. Consider the
regular covering $Y_{G}$ over $Y$ with the transformation group $G$. From the
assumption for $i_{*}$ it follows that
there are a regular covering $X_{G}$ over $X$ with the transformation group
$G$ (the pullback of $Y_{G}$ by $i$) and the
holomorphic embedding $i_{G}:X_{G}\hookrightarrow Y_{G}$ that covers
$i$. Divisor $X\subset Y$ determines a holomorphic line
bundle $L$ over $Y$ and a holomorphic section $s:Y\longrightarrow L$ with
a simple zero along $X$. Further, for every $p\in X$, there is a coordinate
neighbourhood $(U,z)$ centered at $p$ and a holomorphic frame $e$ for $L$  on
$U$ such that $s=z_{1}e$ on $U$. Let $h$ be a hermitian metric on $L$ and
$\nabla$ be the canonical connection with curvature $\Theta$ constructed by
$h$. By the same letters we denote the pullback of $L$, $h$, $s$ and $\Theta$
to $Y_{G}$. Note also that if $\phi$ is a smooth function on $Y_{G}$ then the
weighted metric $e^{\phi}h$ on $L$ has a curvature
$\Theta'=-{\cal L}(\phi)+\Theta$.

Let $U_{0}$ be the pullback of the complement of a closed neighbourhood of
$X\subset Y$ and $U_{1},...,U_{N}$ be the pullbacks of
shrunk coordinate polydisks covering a larger neighbourhood of $X$.
Also pull back a smooth partition of unity $(\xi_{i})$ subordinate to
$(U_{i})$. Let $f$ be a holomorphic function on $X_{G}$ such that
$f^{2}e^{-\phi}$ is integrable on $X_{G}$. For $i\geq 1$, extend $f$ to a
holomorphic function $f_{i}$ on $U_{i}$ which is constant on each line
$\{z_{2},...,z_{n}\ \ {\rm constant}\}$. Let $f_{0}=0$ on $U_{0}$.
Since $f_{i}=f=f_{j}$ on $X_{G}$ and $X_{G}$ is smooth, we
can define a holomorphic section of the dual bundle $L^{*}$ on
$U_{ij}=U_{i}\cap U_{j}$ by the formula
$$
u_{ij}=(f_{i}-f_{j})s^{-1}\ .
$$
Then
$$
v_{i}=\sum_{i}u_{ij}\xi_{j}
$$
is a smooth section of $L^{*}$ on $U_{j}$ and $v_{i}-v_{j}=u_{ij}$.
Hence $\overline{\partial}v_{i}=\overline{\partial}v_{j}$ on $U_{ij}$,
so we get a $\overline{\partial}$-closed, $L^{*}$-valued (0,1)-form
$\eta$ on $Y_{G}$ defined as $\overline{\partial}v_{i}$ on $U_{i}$.
Assume that $\phi$ satisfies (\ref{lip}) or (\ref{ineq}), where $d$
is the distance on $Y_{G}$ defined by the pullback of a metric
on $Y$. Denote by $|f|$ the weighted $L_{2}$-norm of $f$ with the weight
$e^{-\phi}$.
\begin{Lm}\label{eta}
(1)\ If
$\phi$ satisfies (\ref{lip}) then $\eta\in L_{2}^{0,1}(Y_{G},L^{*})$ for
$L$ equipped with the metric $e^{\phi}h$ and
$|\eta|\leq C(X,Y,h,\phi)|f|$ in the corresponding  $L_{2}$-norms.\\
(2)\ If $\phi$ satisfies (\ref{ineq}) then
$\eta\in L_{2}^{0,1}(Y_{G},L^{*})$ for $L$ equipped with the metric
$e^{(1+\epsilon)\phi}h$, $\epsilon>0$, and
$|\eta|\leq C(X,Y,h,\phi,\epsilon)|f|$.
\end{Lm}
{\bf Proof.} We prove (2). The proof of (1) goes along the same
lines (see also arguments in [La, Th. 3.1]).

We have to show that $|\eta|^{2}e^{-(1+\epsilon)\phi}$ is integrable on
$Y_{G}$. On $U_{0}$, $s$ is bounded away from 0 and
$$
\eta=\overline{\partial}v_{0}=-\sum_{j}f_{j}s^{-1}\overline{\partial}\xi_{j}\ ,
$$
so
$$
|\eta|^{2}\leq c\sum_{j}|f_{j}|^{2},
$$
where $c$ depends only on $X,Y,h$. Further,
$$
\int_{U_{j}}|f_{j}|^{2}e^{-(1+\epsilon)\phi}\omega^{n}\leq c'(\epsilon,X,Y,
\phi)\int_{X_{G}\cap U_{j}}|f|^{2}e^{-\phi}\omega^{n-1}
$$
because $\phi$ satisfies (\ref{ineq}).
Since $f^{2}e^{-\phi}$ is integrable on $X_{G}$, so is
$|\eta|^{2}e^{-(1+\epsilon)\phi}$ on $U_{0}$. For $i\geq 1$,
$$
\eta=\overline{\partial} v_{i}=\sum_{j}(f_{i}-f_{j})s^{-1}
\overline{\partial}\xi_{j}
$$
on $U_{i}$ and it remains to show that
\begin{equation}\label{star}
\sum_{i,j\geq 1}\int_{U_{ij}}|f_{i}-f_{j}|^{2}|s|^{-2}e^{-(1+\epsilon)\phi}
\omega^{n}<\infty\ .
\end{equation}
For $x\in U_{ij}$, $i,j\geq 1$, there are $x_{i}\in X_{G}\cap U_{i}$ and
$x_{j}\in X_{G}\cap U_{j}$ such that $f_{i}(x)=f(x_{i})$,
$f_{j}(x)=f(x_{j})$ and $d(x_{i},x_{j})\leq c(h,X,Y)|s(x)|$. So,
$$
|f_{i}(x)-f_{j}(x)||s(x)|^{-1}\leq c'(X,Y,d)\sup |df|,
$$
where supremum is taken over $X_{G}\cap (U_{i}\cup U_{j})$. By the Cauchy
inequalities and since $\phi$ satisfies (\ref{ineq}),
$$
\begin{array}{c}
\displaystyle
\int_{U_{ij}}|f_{i}-f_{j}|^{2}|s|^{-2}e^{-(1+\epsilon)\phi}\omega^{n}\leq
c'(X,Y,d)\int_{U_{ij}}\sup|df|^{2}e^{-(1+\epsilon)\phi}\omega^{n}\\
\displaystyle \leq c''(X,Y,d,\epsilon)
\int_{X_{G}\cap (V_{i}\cup V_{j})}|f|^{2}e^{-\phi}\omega^{n-1},
\end{array}
$$
where $V_{i}\supset U_{i}$, $V_{j}\supset U_{j}$ are pullbacks of larger
polydisks. Since $f^{2}e^{-\phi}$ is integrable on $X_{G}$, (\ref{star})
follows.

The lemma is proved.\ \ \ \ \ $\Box$

Assume now that under conditions of Lemma \ref{eta} there is a smooth section
$w$ of $L^{*}$ such that $\overline{\partial}w=\eta$ and $|w|^{2}e^{-\phi}$
(respectively, $|w|^{2}e^{-(1+\epsilon)\phi}$) is integrable.
Let $u_{i}=v_{i}-w$.
Then $u_{i}$ is a holomorphic section of $L^{*}$ on $U_{i}$ and
$u_{i}-u_{j}=u_{ij}$, so
$$
f_{i}-u_{i}\otimes s=f_{j}-u_{j}\otimes s\ \ \ \ {\rm on}\ \ \ U_{ij}\ .
$$
Hence we obtain a holomorphic extension $F$ of $f$ to $Y$ by setting
$$
F=f_{i}-u_{i}\otimes s=f_{i}+w\otimes s-\sum_{j}(f_{i}-f_{j})\xi_{j}\ \ \
{\rm on}\ \ U_{i}.
$$
The term $w\otimes s$ is $L_{2}$ with respect to $e^{-\phi}$
(respectively, $e^{-(1+\epsilon)\phi}$) by construction of $w$ and since $s$
is bounded.
The other two terms on the right-hand side can be shown to be $L_{2}$ with
respect to $e^{-\phi}$ (respectively, $e^{-(1+\epsilon)\phi}$) by arguments
similar to those used for $\eta$ above. Hence $F^{2}e^{-\phi}$ (respectively,
$F^{2}e^{-(1+\epsilon)\phi}$) is integrable.\\
{\bf 2.3. Symmetric products of curves.}
Let $\Gamma$ be a complex compact curve of genus $g\geq 1$,
$\Gamma^{\times g}$ and $S\Gamma^{\times g}$ be the direct and the
symmetric products of $g$-copies of $\Gamma$. Then the manifold
$S\Gamma^{\times g}$ is the quotient of $\Gamma^{\times g}$ by the action of
the permutation group $S_{g}$. Therefore there exists a finite holomorphic
surjective map $\pi:\Gamma^{\times g}\longrightarrow S\Gamma^{\times g}$.
Further, $S\Gamma^{\times g}$ is birational isomorphic to $\Co\To^{g}$
(denote this isomorphism by $j$). Let $(p,...,p)\in\Gamma^{\times g}$
be a fixed point. Denote by $\Gamma^{k}$, $k\leq g$,
submanifold
$\{(p,...,p, z_{1},...,z_{k})|\ z_{1},...,z_{k}\in\Gamma\}
\subset\Gamma^{\times g}$.
\begin{Lm}\label{sym}
For any $k$, image $\pi(\Gamma^{k})$ is a complex submanifold of
$S\Gamma^{\times g}$.
\end{Lm}
{\bf Proof.}
For a point $y=(p,...,p,z_{1},...,z_{k})\in\Gamma^{k}$ consider its orbit
$o(y):=S_{g}(y)$. By definition, $\pi$ maps $o(y)$ to $\pi(y)$ and
intersection $o(y)\cap\Gamma^{k}=\{(p,...,p,S_{k}(z))\}$; here
$z=(z_{1},...,z_{k})$ and $S_{k}$ is the permutation group acting on the
set of $k$ elements. The quotient by the action of $S_{k}$ is
manifold $X_{k}:=(p,...,p,S\Gamma^{\times k})$. So we have a holomorphic
injective mapping $\pi_{k}:X_{k}\longrightarrow S\Gamma^{\times g}$ whose
image coincides with $\pi(\Gamma^{k})$. Now let $y_{0}$ be local
coordinates in a neighbourhood $U_{0}$ of
$p\in\Gamma$ and $y_{i}$, $1\leq i\leq k$,
be local coordinates in a neighbourhood $U_{i}$ of $z_{i}\in\Gamma$ such that
$y_{0}(p)=0, \ y_{i}(U_{i})\cap y_{j}(U_{j})=\emptyset \ {\rm for}\
z_{i}\neq z_{j} \ {\rm and}\ y_{i}=y_{j} \ {\rm in}\ U_{i}=U_{j}\
{\rm for}\ z_{i}=z_{j}$. Denote by $\sigma_{1},...,\sigma_{g}$ elementary
symmetric functions from $g$ variables. For $(z_{1},...,z_{g})\in
U_{0}\times...\times U_{0}\times U_{1}\times...\times U_{k}\subset
\Gamma^{\times g}$ set
$u_{i}(z)=y_{0}(z_{i})$, $1\leq i\leq g-k$, and $u_{i}(z)=y_{i}(z_{g-k+i})$,
$1\leq i\leq k$.
By the theorem on symmetric polynomials the mapping
$$
f:(w_{1},...,w_{g})\mapsto (\sigma_{1}(u(w)),...,\sigma_{g}(u(w)))
$$
determines a local coordinate system on
$\pi(U_{0}\times...\times U_{0}\times U_{1}\times...\times U_{k})\subset
S\Gamma^{\times g}$ (see [GH, Ch. 2, p. 259]). Then the image of
restriction $f|_{\pi(\Gamma_{k})}$ belongs to $\Co^{k}\subset\Co^{g}$.
By the same reason $f\circ\pi_{k}$ determines a local coordinate
system in the corresponding neighbourhood on $X_{k}$. This shows that
$\pi_{k}$ is a biholomorphic embedding. Thus we proved that $\pi(\Gamma_{k})$
is smooth.\ \ \ \ \ $\Box$\\
{\bf 2.4. Norm estimates.} Let $M$ and $N$ be compact Riemannian manifolds
and $f:M\longrightarrow N$ be a smooth surjective map. Assume that
$f_{*}:\pi_{1}(M)\longrightarrow\pi_{1}(N)$ is a surjection. Let $G$
be a quotient group of $\pi_{1}(N)$ and $N_{G}$, $M_{G}$ regular coverings
with the transformation group $G$ over $N$ and $M$, respectively, such that
$M_{G}=f^{*}N_{G}$. Then there is a map
$f_{G}:M_{G}\longrightarrow N_{G}$ that covers $f$. We consider $M_{G}$ and
$N_{G}$ in the metrics pulled back from $M$ and $N$, respectively.
Further, if $E^{p}(K)$ is the space of $p$-forms on a Riemannian manifold $K$
denote by $|\cdot|_{x}$ the norm in the vector space $E^{p}(K)|_{x} (\cong
\wedge^{p} T_{x}^{*})$, $x\in K$, constructed by the metric dual to the
Riemannian one.
\begin{Lm}\label{compare}
Let $\omega$ be a bounded differential $p$-form on $N_{G}$, i.e.,
$\sup_{x\in N_{G}}|\omega|_{x}<\infty$.
Then there is $C=C(f,p)>0$ such that
$$
|f_{G}^{*}(\omega)|_{x}\leq C|\omega|_{f_{G}(x)}\ .
$$
\end{Lm}
{\bf Proof.} Let us write $\omega$ in local orthogonal coordinates
lifted from $N$.  Then the compactness arguments show that the
statement follows easily from a similar statement for elements of
the orthogonal basis. We leave the details to the reader. \ \ \ \
\ $\Box$
%=================================
\sect{\hspace*{-1em}. Proofs.}
We prove Theorem \ref{te2} only. The proof of Theorem \ref{te1}
is similar and can be obtained by removing $\epsilon$ in the arguments below.\\
{\bf 3.1.} We start by proving Theorems \ref{te1} and \ref{te2} for
curves. \\
{\bf Proof of Theorem \ref{te2} for curves.}
Assume that the Albanese map $A:\Gamma\longrightarrow\Co\To^{g}$ is defined
with respect to a basic point $p\in\Gamma$.
For $X_{i}:=\pi(\Gamma^{i})\subset S\Gamma^{\times g}$ consider
the flag of submanifolds $X_{1}\subset...\subset X_{g}=S\Gamma^{\times g}$
(see definitions in Section 2.3). The Jacobi map
$j:S\Gamma^{\times g}\longrightarrow\Co\To^{g}$ maps, by definition,
$X_{1}$ biholomorphically to $A(\Gamma)$
(which we identify with $\Gamma$). Moreover, the fundamental group
$\pi_{1}(S\Gamma^{\times g})$ is isomorphic (under $j_{*}$) to
$\pi_{1}(\Co\To^{g})=\Z^{2g}$ and embedding
$X_{i}\subset S\Gamma^{\times g}$ induces a surjective homomorphism
of fundamental groups. Thus if $G$ is a quotient group of
$\pi_{1}(\Co\To^{g})$
one can construct regular coverings $X_{iG}$ over
$X_{i}$, $i=1,...,g$, with transformation group $G$ such that
$X_{1G}\subset....\subset X_{gG}$ is a flag of complex submanifolds covering
the flag $X_{1}\subset...\subset X_{g}$ and
there is a proper surjective map with connected fibres
$j_{G}:X_{gG}\longrightarrow T_{G}$ that covers $j$.

For any function $f\in {\cal H}_{\phi}(\Gamma_{G})$ consider its
pullback $f_{1}:=j_{G}^{*}(f)$ on $X_{1G}$. Then according to
Lemma \ref{compare}, $f_{1}$ belongs to the space ${\cal
H}_{j_{G}^{*}(\phi)}(X_{1G})$ determined with respect to the
pullback of the volume form of $X_{1}$. Moreover,
$j_{G}^{*}(\phi)$ satisfies condition (\ref{ineq}) (respectively,
(\ref{lip})) for the distance $d'$ defined by the pullback of a
K\"{a}hler metric on $X_{g}$. It follows from the inequality $$
d(j_{G}(x),j_{G}(y))\leq C(j_{G})d'(x,y)\ \ \ (x,y\in X_{gG}). $$
Now for a sufficiently small $\epsilon>0$ we prove that $f_{1}$
admits an extension $f_{2}\in {\cal
H}_{(1+\epsilon)j_{G}^{*}(\phi)}(X_{2G})$ satisfying conditions of
Theorem \ref{te2}; $f_{2}$  admits a similar extension $f_{3}\in
{\cal H}_{(1+2\epsilon)j_{G}^{*}(\phi)}(X_{3G})$ etc. Finally, we
obtain an extension $f_{g}\in {\cal
H}_{(1+(g-1)\epsilon)j_{G}^{*}(\phi)}(X_{gG})$ of $f_{1}$.
Clearly, $f_{g}$ is constant on fibres of $j_{G}$ and thus
determines a function $f'\in {\cal
H}_{(1+(g-1)\epsilon)\phi}(T_{G})$ that extends $f$. Our arguments
will guarantee its uniqueness and fulfillment of the required norm
estimates. This will finish the proof.

We use inductive arguments. Assume that we have the required extension
$f_{k}\in {\cal H}_{(1+(k-1)\epsilon)j_{G}^{*}(\phi)}(X_{kG})$ of $f_{k-1}$.
Construct now extension
$f_{k+1}\in{\cal H}_{(1+k\epsilon)j_{G}^{*}(\phi)}(X_{(k+1)G})$.

For each $k$ consider the regular covering $Y_{k}$ over
$\Gamma^{k}\subset\Gamma^{\times g}$ with the transformation group
$G$. Since the map $j\circ\pi:\Gamma^{k}\longrightarrow\Co\To^{g}$
is invariant with respect to the action of the permutation group
$S_{k}$ acting on $\Gamma^{k}(\cong\Gamma^{\times k})$ and
$(p,...,p)\in\Gamma^{\times g}$ is a fixed point with respect to
$S_{k}$, by the covering homotopy theorem there is a covering
action of $S_{k}$ on $Y_{k}$. Moreover, there is a holomorphic map
$\pi_{G}:Y_{g}\longrightarrow X_{gG}$ that covers
$\pi:\Gamma^{\times g}\longrightarrow S\Gamma^{\times g}$ and is
invariant with respect to the action of $S_{g}$. Consider the
orbit $V_{k}=S_{k+1}(Y_{k})$ in $Y_{k+1}$. Then $V_{k}$ covers the
orbit $W_{k}=S_{k+1}(\Gamma^{k})\subset\Gamma^{k+1}$.
\begin{Lm}\label{plus}
Divisor $W_{k}$ determines a positive line bundle $E_{k}$ over
$\Gamma^{k+1}$.
\end{Lm}
{\bf Proof.} Assume without loss of generality that $\Gamma^{k+1}=
\Gamma^{\times (k+1)}$. Let $P:\Gamma^{k+1}\longrightarrow\Gamma$
be the projection defined by $$ P(z_{1},...,z_{k+1})=z_{1},\ \ \
\\ (z_{1},...,z_{k+1})\in\Gamma^{\times (k+1)}\ . $$
Then $P^{-1}(x)=(x,\Gamma^{\times k})$ for a fixed $x\in\Gamma$.
Denote by $E_{x}$ a positive line bundle over $\Gamma$ defined by
the divisor $\{x\}$ and by $\Theta_{x}$ its curvature (for a
suitable hermitian metric on $E_{x}$) such that
$\frac{\sqrt{-1}}{2\pi}\Theta_{x}$ is a positive (1,1)-form. Let
$e_{i}\in S_{k+1}$, $i=1,...,k+1$, be such that
$\cup_{i}e_{i}^{-1}(\Gamma^{k})=S_{k+1}(\Gamma^{k})$. Then by
definition, $E_{k}=\otimes_{i}e_{i}^{*}(P^{*}E_{x})$ is a positive
line bundle over $\Gamma_{k+1}$. In fact, if in local coordinates
$P^{*}\Theta_{x}=a(z_{1},\overline{z_{1}})dz_{1}\wedge
d\overline{z}_{1}$ with $a(z_{1},\overline{z_{1}})>0$, the
curvature $\Theta_{k}$ of $E_{k}$ equals
$\sum_{i=1}^{k}a(z_{i},\overline{z_{i}})dz_{i}\wedge
d\overline{z}_{i}$. Clearly, $\frac{\sqrt{-1}}{2\pi}\Theta_{k}$ is
positive implying that $E_{k}$ is positive.\ \ \ \ \ $\Box$

Let $h_{k}$ be a hermitian metric on $E_{k}$ with the curvature
$\Theta_{k}$. By the same letters we denote the pullback of
$h_{k}$, $E_{k}$ and $\Theta_{k}$ to $Y_{k+1}$. Let $L_{k}$ be the
holomorphic vector bundle on $X_{k+1}$ defined by the divisor
$X_{k}$ and $h_{k}'$ a hermitian metric on $L_{k}$. By the same
letters  we also denote the pullback of $L_{k}$ and $h_{k}'$ to
$X_{(k+1)G}$.  Below we consider $L_{k}$ with the weighted metric
$e^{(1+k\epsilon)j_{G}^{*}(\phi)}h_{k}'$. By Lemma \ref{eta} (2)
there is a linear continuous mapping $F_{k,\epsilon}:{\cal
H}_{(1+(k-1)\epsilon)j_{G}^{*}(\phi)}(X_{kG}) \longrightarrow
Z_{2}^{0,1}(X_{(k+1)G},L_{k}^{*})$. Put
$\eta_{k}=F_{k,\epsilon}(f_{k})$. Since, by definition,
$\pi^{-1}(X_{k})\cap\Gamma^{k+1}=W_{k}$, the bundle
$\pi_{G}^{*}L_{k}$ equals $E_{k}$ on $\Gamma^{k+1}$. In
particular, $\eta_{k}'=\pi_{G}^{*}(\eta_{k})$ is a
$\overline{\partial}$-closed (0,1)-form on $Y_{k+1}$ with values
in $E_{k}^{*}$ and $\eta_{k}'\in L_{2}^{0,1}(Y_{k+1},E_{k}^{*})$
for $E$ equipped with the metric $e^{(1+k\epsilon)\phi'}h_{k}$
where $\phi'=\pi_{G}^{*}(j_{G}^{*}(\phi))$. Further, the curvature
${\cal R}_{k}$ of $E_{k}$ equals $-(1+k\epsilon){\cal
L}(\phi')+\Theta_{k}$. Moreover, according to Lemma \ref{compare},
$$ |{\cal L}(\phi')|_{x}\leq C(j_{G}\circ\pi_{G},2)||{\cal
L}(\phi)|_{(j_{G}\circ\pi_{G})(x)}\ \ \ \ (x\in Y_{g}).  $$ In
particular, there is a positive constant $C$ (depending on
$\Gamma$ only) such that

{\em if $\sup_{x\in T_{G}}|{\cal L}(\phi)|_{x}<C'<C$ and
$0<\epsilon\leq 1/g$, $1\leq k\leq g$, there is an $a=a(C')>0$ so
that ${\cal R}_{k}>a\Theta_{k}$.}\\ Let $\phi$ satisfy the above
condition and $\epsilon<1/g$. Since $\Theta_{k}$ is a K\"{a}hler
form on $\Gamma_{k+1}$, according to Proposition \ref{lar} and
Remark \ref{lin} there is a linear continuous mapping
$s_{k}:Z_{2}^{0,1}(Y_{k+1},E_{k}^{*})\longrightarrow
W_{2}^{0,0}(Y_{k+1},E_{k}^{*})$ inverse to $\overline{\partial}$.
Then for $r_{k}=s_{k}(\eta_{k}')$ we have
$\overline{\partial}r_{k}=\eta_{k}'$ and $r_{k}\in
L_{2}(Y_{k+1},E_{k}^{*})$. Applying now arguments similar to those
used in Section 2.2 (for the pullback to $Y_{k+1}$ of local
extensions of $f_{k}$) get a holomorphic function $g_{k+1}$ on
$Y_{k+1}$ that extends $\pi_{G}^{*}(f_{k})$ and belongs to ${\cal
H}_{(1+k\epsilon)\phi'}(Y_{k+1})$.

Assume also that there is another extension $g'\in {\cal
H}_{(1+k\epsilon)\phi'}(Y_{k+1})$ of $\pi_{G}^{*}(f_{k})$. Let
$s'$ be the pullback to $Y_{k+1}$ of a holomorphic section of the
bundle $E_{k}$ on $\Gamma_{k+1}$ with a simple zero along $W_{k}$.
(Recall that the pullback of $E_{k}$ we denote by the same
letter). Then $d=(g_{k+1}-g')(s')^{-1}$ is an $L_{2}$ integrable
holomorphic section of $E_{k}^{*}$. Here $E_{k}$ is taken with the
weighted metric $e^{(1+k\epsilon')\phi'}h_{k}$, where an
$\epsilon'$ satisfies $\epsilon<\epsilon'<1/g$. The arguments are
similar to those used in the proof of Lemma \ref{eta}.  Therefore
according to Proposition \ref{lar} for $q=0$, the function $d$ is
zero. This proves the uniqueness of the extension. Since
$\pi_{G}^{*}(f_{k})$ is invariant with respect to the action of
the permutation group $S_{k}$, for any $e\in S_{k}$ the function
$e^{*}(g_{k+1})$ is also an extension of $\pi_{G}^{*}(f_{k})$
belonging to ${\cal H}_{(1+k\epsilon)\phi'}(Y_{k+1})$. Thus the
uniqueness of extension implies that $e^{*}(g_{k+1})=g_{k+1}$.  So
there is a uniquely defined holomorphic function $f_{k+1}$ on
$X_{(k+1)G}$ such that $\pi_{G}^{*}(f_{k+1})=g_{k+1}$, $f_{k+1}\in
{\cal H}_{(1+k\epsilon)j_{G}^{*}(\phi)}(X_{k+1}G)$ and $f_{k+1}$
is an extension of $f_{k}$. In fact our arguments (based on Remark
\ref{lin}) show that we constructed a linear continuous extension
operator which gives us the required norm estimates. Therefore, by
induction, we get a holomorphic function $f_{g}$ on $X_{gG}$ which
belongs to ${\cal H}_{(1+(g-1)\epsilon)j_{G}^{*}(\phi)}(X_{gG})$
and extends $f_{1}$. As it was noted at the beginning of the
proof, $f_{g}$ determines the required extension of $f$. This
proves Theorem \ref{te2} for curves.\ \ \ \ \ $\Box$\\ {\bf Proof
of Theorem \ref{te2} for projective manifolds.} Let $M$ be a
projective manifold of dimension $n\geq 2$ with a very ample line
bundle $L$ and with a K\"{a}hler form $\omega$. We may think of
$M$ as embedded in some projective space and of $L$ as the
restriction to $M$ of the hyperplane bundle with the standard
positively curved metric. Then zero loci of sections of $L$ are
hyperplane sections of $M$. By Bertini's theorem, the generic
linear subspace of codimension $n-1$ intersects $M$ transversely
in a $\penalty-10000$ smooth curve $C$. By the Lefschetz
hyperplane theorem, $C$ is connected and the map
$\pi_{1}(C)\longrightarrow\pi_{1}(M)$ is surjective. Let $M_{G}$
be the regular covering over $M$ with a free abelian
transformation group $G$. Then the regular covering $C_{G}$ over
$C$ with the same transformation group $G$ is embedded into
$M_{G}$. Assume that $f\in {\cal H}_{\tilde\phi}(M_{G})$ with
$\tilde\phi$ satisfying (\ref{ineq}). Then $g:=f|_{C_{G}}$ belongs
to ${\cal H}_{(1+\epsilon)\tilde\phi}(C_{G})$ for any positive
$\epsilon$. Indeed, let $U_{1},...,U_{N}$ be the pullbacks to
$M_{G}$ of shrunk coordinate polydisks covering an open
neighbourhood of $C\subset M$ and $V_{i}\supset U_{i}$ be
pullbacks of larger polydisks.  We may assume that $C_{G}\cap
V_{i}=\{z_{1}=0\}$, $i=1,...,N$, for the pullback of the
corresponding local coordinates. Then application of (\ref{ineq})
and subharmonicity of $|f|^{2}$ get $$ \int_{C_{G}\cap
U_{i}}|f|^{2}e^{-(1+\epsilon)\tilde\phi}\omega\leq
c(M)\int_{V_{i}}|f|^{2}e^{-\tilde\phi}\omega^{n}<\infty\ . $$ This
implies $g\in {\cal H}_{(1+\epsilon)\tilde\phi}(C_{G})$. Let
$C=M_{1}\subset M_{2}\subset ...\subset M_{n}=M$ be a flag of
projective submanifolds of $M$, where $M_{i}$ is intersection of
$M$ with the generic linear subspace of codimension $n-i$. Let
$C_{G}=M_{1G}\subset...\subset M_{nG}=M_{G}$ be the flag of the
corresponding regular coverings with the transformation group $G$.
Then the arguments similar to those used in Section 3.1 (see also
arguments in Theorem 3.1 of [La]) show that if $L$ is very ample
then $g$ admits a unique extension $f'\in {\cal
H}_{(1+\epsilon+\delta)\tilde\phi}(M_{G})$ for a sufficiently
small positive $\epsilon$ and $\delta=\delta(\epsilon)$. But
clearly in this case $f=f'$. Thus we proved that $f$ is uniquely
determined by $f|_{C_{G}}$.

Let now $A:M\longrightarrow\Co\To^{k}$ be the Albanese map for $M$
defined with respect to a point $p\in C$ by integration of
holomorphic 1-forms $\omega_{1},...,\omega_{k}\in\Omega^{1}(M)$
(generating a basis there). Set $\eta_{i}:=\omega_{i}|_{C}$ for
$i=1,...,k$. Then by the Lefschetz theorem $\eta_{1},...,\eta_{k}$
are linearly independent in $\Omega^{1}(C)$. Choose 1-forms
$\eta_{k+1},...,\eta_{s}\in \Omega^{1}(C)$ such that
$\eta_{1},...,\eta_{s}$ generates a basis. Further, define the
Albanese map $A':C\longrightarrow\Co\To^{s}$ with respect to the
point $p$ by integration the forms of this basis. Then according
to our construction there is a surjective map
$P:\Co\To^{s}\longrightarrow\Co\To^{k}$ whose fibres are complex
tori such that $P_{*}:\pi_{1}(\Co\To^{s})\longrightarrow
\pi_{1}(\Co\To^{k})$ is a surjection and $A=P\circ A'$. Denote by
$T_{G}'$ a regular covering over $\Co\To^{s}$ with the
transformation group $G$. Then there is a complex map
$P_{G}:T_{G}'\longrightarrow T_{G}$ that covers $P$ whose fibres
are also tori. Let $A_{G}':C_{G}\longrightarrow T_{G}'$ be the map
covering $A'$ and $\phi'=P_{G}^{*}(\phi)$. Note that  $\phi'$
satisfies (\ref{ineq}) on $T_{G}'$ and
$(A_{G}')^{*}(\phi')=\tilde\phi|_{C_{G}}$. Applying Theorem
\ref{te2} for curves to the map $A_{G}':C_{G}\longrightarrow
T_{G}'$ and the function $\phi'$ we obtain

{\em there is $C=C(M,A')>0$ such that for $|{\cal L}(\phi')|<C'<C$
and for sufficiently small positive numbers
$\epsilon\leq\epsilon(C')$, $\delta\leq\delta(C')$ there is a
uniquely defined holomorphic function $\tilde f\in {\cal
H}_{(1+\epsilon+\delta)\phi'}(T_{G}')$ satisfying
$g=(A_{G}')^{*}(\tilde f)$ and $|\tilde f|\leq
C(\epsilon,\delta)|g|$ (in the corresponding $L_{2}$-norms).}\\
Since $P_{G}$ is a proper map with connected fibres, $\tilde f$
determines a function $h\in {\cal
H}_{(1+\epsilon+\delta)\phi}(T_{G})$ such that
$A_{G}^{*}(h)|_{C_{G}}=g$ and $|h|\leq \tilde
C(\epsilon,\delta)|g|$. But as we proved, $f$ is uniquely
determined by $g=f|_{C_{G}}$ and $|g|\leq c(\epsilon)|f|$.
Therefore $A_{G}^{*}h=f$ and $h$ satisfies the required norm
estimate. Finally, by Lemma \ref{compare}, $|{\cal L}(\phi')|\leq
c(P)|{\cal L}(\phi)|$ and so the above extension theorem is valid
for any $B'$ satisfying $|{\cal L}(\phi)|<B'<C/c(P)$.

This completes the proof of Theorem \ref{te2} for projective
manifolds. \ \ \ \ \ $\Box$\\ {\bf 3.2.} {\bf Proof of Theorem
\ref{te3}.} Let $M_{G}$ be a regular covering over $M$ with the
transformation group $G$ and $A_{G}:M_{G}\longrightarrow T_{G}$ be
the covering of the Albanese map $A:M\longrightarrow\Co\To^{n}$.
Assume that $f$ is a holomorphic function on $M_{G}$ satisfying $$
|f(z)|\leq c(\epsilon)e^{\epsilon r(z)}\ \ \ (z\in M_{G}) $$ for
any $\epsilon>0$. Let $\phi$ be the distance from a fixed point in
$T_{G}$ in the flat metric pulled back from $\Co\To^{n}$ and
$\tilde\phi=A_{G}^{*}(\phi)$. Further, by $\rho_{G}$ denote the
distance from $0$ on $G(\cong\Z^{k})$ determined with respect to
the word metric. Since by our construction growth of $r$ and
$\tilde\phi$ is equivalent to growth of $\rho_{G}$, the function
$f$ belongs to ${\cal H}_{\epsilon\tilde\phi}(M_{G})$ for any
$\epsilon>0$. We now apply Theorem \ref{te1}. Here we assume that
$|{\cal L}(\phi)|$ is sufficiently small replacing if necessary
$\phi$ by a smooth function $\phi_{1}$ with the same growth such
that $|{\cal L}(\phi_{1})|$ is small. In fact  $\phi_{1}$ can be
constructed as follows.\\ Note, first, that $T_{G}$ is
diffeomorphic to $\To^{2n-k}\times\Re^{k}$ where second
derivatives of the diffeomorphism are bounded in the flat
coordinate system on $T_{G}$. Then put
$\phi_{1}(v,x):=\sqrt{p+|x|^{2}}$, for
$(v,x)\in\To^{2n-k}\times\Re^{k}$, where $|x|$ is the Euclidean
norm of $x\in\Re^{k}$ and $p$ is sufficiently big positive number.

Further, according to Theorem \ref{te1} there is a
uniquely defined holomorphic function $f'\in \cap_{\epsilon>0}
{\cal H}_{\epsilon\phi}(T_{G})$ such that $A_{G}^{*}(f')=f$.
Prove now that $f'$ is a constant.

We regard the maximal free abelian subgroup $\Gamma\subset H_{1}(M,\Z)$ as a
lattice in $\Co^{n}$ determining $\Co\To^{n}$ and $H\subset\Gamma$ as a
sublattice such that the minimal complex vector space containing $H$ is
$\Co^{n}$. Consider the pullback
$g$ of $f'$ to $\Co^{n}$. Clearly $g$ is invariant  with respect to the action
(by shifts) of $H$ and satisfies
$$
|f(z)|\leq c(\epsilon)e^{\epsilon |z|}
$$
for any positive $\epsilon$. For an element $e_{1}\in H$
let $X_{1}$ be a minimal complex vector space containing
$\{ne_{1}\}_{n\in\Z}$. For any $z\in\Co^{n}$ consider restriction
$g'=g|_{z+X_{1}}$. We identify $z+X_{1}$ with $\Co$ and $\{ne_{1}\}_{n\in\Z}$
with $\Z$. Then $g'$ is a
holomorphic function on $\Co$ of an arbitrary small exponential type which is
constant on $\Z$. Therefore by Cawrtright's theorem [Ca],
$g'$ is constant on $\Co$. This implies that $g(z+v)=g(z)$ for any
$z\in\Co^{n}$ and $v\in X_{1}$.  In particular, there is a holomorphic
function $g_{1}$ on the quotient $\Co^{n}/X_{1}=\Co^{n-1}$ of an arbitrary
small exponential type whose pullback to $\Co^{n}$ coincides with $g$.
Denote by $H_{1}$ image of $H$ in $\Co^{n}/X_{1}=\Co^{n-1}$. By definition,
$g_{1}$ is
invariant with respect to the action of $H_{1}$ and the minimal complex vector
space containing $H_{1}$ is $\Co^{n-1}$. Choose
$e_{2}\in H_{1}$ and denote by $X_{2}$ the minimal complex subspace
containing $\{ne_{2}\}_{n\in\Z}$. Applying the very same arguments we get
$g_{1}(z+v)=g_{1}(z)$ for any $z\in\Co^{n-1}$ and $v\in X_{2}$. Continuing
by induction we finally obtain that the initial function $g$ is constant.

This completes the proof of the theorem. \ \ \ \ \ $\Box$

Note that our arguments give a more general statement.
\begin{Th}\label{te4}
Let $H$ be such that $\Lambda (H)=\Omega^{1}(M)$. Then there is a
positive constant $\sigma=\sigma(M)$ such that any holomorphic on
$M_{G}$ function $f$ satisfying $$ |f(z)|\leq ce^{\sigma' r(z)}\ \
\ \ (0<\sigma'<\sigma,\ z\in M_{G})$$ is a constant.
\end{Th}
%===================================================

{\em Department of Mathematics, Ben Gurion University of the
Negev, P.O.Box 653,

Beer-Sheva 84105, Israel}

{\em E-mail address: brudnyi@cs.bgu.ac.il
\end{document}
%==================================
%==================================
%==================================